\documentclass[12pt,leqno]{article}
\usepackage{amsmath, amsthm, amsfonts, amssymb, color}
\usepackage{mathrsfs}
\theoremstyle{definition}

\newcommand{\scr}[1]{\mathscr #1}
\definecolor{wco}{rgb}{0.5,0.2,0.3}

\numberwithin{equation}{section} \theoremstyle{remark}

\newcommand{\ua}{\uparrow}

\title{\bf   Intrinsic Ultracontractivity on Riemannian Manifolds with Infinite Volume Measures}
\author{
  {\bf Feng-Yu Wang}\footnote{Supported in part by Creative Research Group Fund of the National
 Foundation of China (No. 10121101), RFDP(20040027009) and the 973-project in China.}\\
\footnotesize{School  of Mathematical Science \&  Lab. Math. Com. Sys.,}\\ 
\footnotesize{ Beijing Normal
University, Beijing 100875, China}\\
 \footnotesize{ Email: wangfy@bnu.edu.cn}\\
{\footnotesize Current address: WIMCS, Department of Mathematics,
University of Wales Swansea,} \\
{\footnotesize Singleton Park, Swansea SA2 8PP, UK. }}

\begin{document}
\maketitle
\def\R{\mathbb R} \def\Z{\mathbb Z} \def\ff{\frac} \def\ss{\sqrt}
\def\N{\mathbb N}
\def\dd{\delta} \def\DD{\Delta} \def\vv{\varepsilon} \def\rr{\rho}
\def\<{\langle} \def\>{\rangle} \def\GG{\Gamma} \def\ggg{\gamma}
\def\ll{\lambda} \def\LL{\Lambda} \def\nn{\nabla} \def\pp{\partial}
\def\d{\text{\rm{d}}} \def\bb{\beta} \def\aa{\alpha} \def\D{\scr D}
\def\E{\scr E} \def\si{\sigma} \def\ess{\text{\rm{ess}}} \def\Sec{\text{\rm{Sec}}}
\def\beg{\begin} \def\beq{\begin{equation}}  \def\F{\scr F}
\def\Ric{\text{\rm{Ric}}} \def\Hess{\text{\rm{Hess}}}\def\B{\scr B}
\def\e{\text{\rm{e}}} \def\ua{\underline a} \def\OO{\Omega} \def\b{\mathbf b}
\def\oo{\omega}     \def\tt{\tilde} \def\Ric{\text{\rm{Ric}}}
\def\cut{\text{\rm{cut}}} \def\P{\mathbb P} \def\ifn{I_n(f^{\bigotimes n})}
\def\fff{f(x_1)\dots f(x_n)} \def\ifm{I_m(g^{\bigotimes m})} \def\ee{\varepsilon}
\def\pm{\pi_{\mu}}   \def\p{\mathbf{p}}   \def\ml{\mathbf{L}}
 \def\C{\scr C}      \def\aaa{\mathbf{r}}     \def\r{r}
\def\gap{\text{\rm{gap}}} \def\prr{\pi_{\mu,\varrho}}  \def\r{\mathbf r}
\def\Z{\mathbb Z} \def\vrr{\varrho} \def\Sect{{\rm Sect}}
\def\ii{{\rm i}_{\pp M}} \def\Ent{{\rm Ent}} \def\EE{\mathbb E}
\begin{abstract} By establishing the intrinsic super-Poincar\'e
inequality, some explicit conditions are presented for diffusion
semigroups on a non-compact complete Riemannian manifold to be
intrinsically ultracontractive. These conditions, as well as the
resulting uniform upper bounds on the intrinsic heat kernels, are
sharp for some concrete examples.
\end{abstract}

\

 \noindent
 AMS subject Classification:\ \ 58G32, 60J60\\

\noindent
 Keywords: Intrinsic ultracontractivity,  intrinsic super-Poincar\'e inequality, Riemannian manifold, diffusion semigroup.
 \vskip 2cm

\section{Introduction}

Let $(E,\scr F,\mu)$ be a $\si$-finite measure space, and
$(L,\D(L))$ a negative self-adjoint operator generating a
(sub-)Markov semigroup $P_t:= \e^{tL}$ on $L^2(\mu)$. According to
\cite{DS}, the semigroup $P_t$ is called ultracontractive if
$\|P_t\|_{L^1(\mu)\to L^\infty(\mu)}<\infty$ for any $t>0.$ Due to
Gross \cite{G}, the ultracontractivity of $P_t$ can be derived from
log-Sobolev inequalities and for $\mu$ finite,  implies that $P_t$
is compact and  $L$ has empty essential spectrum (see \cite{W00} for
more general results on functional inequalities and the essential
spectrum). On the other hand, however, when $\mu$ is infinite, there
is no direct relationship between the ultracontractivity and  the
spectrum of $L$. For instance, the heat semigroup on $\R^d$ is
ultracontractive with respect to the Lebesgue measure, but the
spectrum of $-\DD$ is continuous. This means that the
ultracontractivity is no longer $``$intrinsic$"$ for the spectrum
property when $\mu$ is infinite. For this reason and for other
applications in the study of Markov semigroups, the intrinsic
ultracontractivity was introduced (cf. \cite{DS}).

Assume that $\ll_0:=\inf\si(-L)$ is a simple eigenvalue of $-L$ with
$\varphi_0>0$ the unique unit eigenfunction, where and in what
follows, we use $\si(\cdot)$ and $\si_{ess}(\cdot)$ to denote the
spectrum and the essential spectrum of an operator respectively. In
general, for the case that $\ll_0$ is an eigenvalue, we  may always
take a nonnegative eigenfunction $\varphi_0$ (called the ground
state in the literature). Indeed, let $\varphi$ be an eigenfunction
with respect to $\ll_0$ with $\mu(\varphi>0)>0$ (otherwise, use
$-\varphi$ in stead of $\varphi$), we have

$$P_t\varphi^+\ge (P_t\varphi)^+=\e^{-\ll_0 t}\varphi^+$$
 for all
$t\ge 0$. But due to the definition of $\ll_0$, one has

$$\|P_t\varphi^+\|_{L^2(\mu)}\le \e^{-\ll_0
t}\|\varphi^+\|_{L^2(\mu)},$$ then $P_t\varphi^+ =\e^{-\ll_0
t}\varphi^+$ for any $t\ge 0$, so that $\varphi^+$ is an
eigenfunction with respect to $\ll_0$ too. If in addition that $P_t$
is irreducible in the sense that $\mu(1_AP_t 1_B)>0$ for any $A,B\in
\F$ with $\mu(A)\mu(B)>0$, then the nonnegative eigenfunction has to
be strictly positive  and $\ll_0$ is a simple eigenvalue. In this
case, $\mu_{\varphi_0}:= \varphi_0^2\mu$ is a probability measure
and

$$P_t^{\varphi_0}f:= \ff {\e^{\ll_0 t}}{\varphi_0} P_t
(f\varphi_0),\ \ \ t\ge 0, f\in L^2(\mu_{\varphi_0})$$ gives rise to
a symmetric $C_0$ Markov semigroup on $L^2(\mu_{\varphi_0}).$ We
call $P_t$ intrinsically ultracontractive if

$$\|P_t^{\varphi_0}\|_{L^1(\mu_{\varphi_0})\to
L^\infty(\mu_{\varphi_0})} <\infty,\ \ \ t>0.$$ Moreover, $P_t$ is
called intrinsically hypercontractive if there exists $t>0$ such
that

$$\|P_t^{\varphi_0}\|_{L^2(\mu_{\varphi_0})\to
L^4(\mu_{\varphi_0})} =1.$$

The intrinsic ultracontractivity has been well studied in the
framework of Dirichlet heat semigroups on (in particular, bounded)
domains in $\R^d$. For instance, the Dirichlet heat semigroup on a
bounded H\"older domain of order $0$ is intrinsically
ultracontractive (see \cite{C, CS}). See the recent work \cite{OW}
and references within for sharp estimates on
$\|P_t^{\varphi_0}\|_{L^1(\mu_{\varphi_0})\to
L^\infty(\mu_{\varphi_0})}$, and   \cite{KS} and references within
for the study of the intrinsic ultracontractivity of L\'evy (in
particular, stable) processes  on domains.

In this paper, we aim to study the intrinsic ultracontractivity for
diffusion semigroups on Riemannian manifolds with infinite invariant
measures. Let $M$ be a complete, connected, non-compact Riemannian
manifold of dimension $d$. Let $L= \DD+\nn V$ for some $V\in
C^2(M)$. Then $L$ generates a unique (Dirichlet) diffusion semigroup
$P_t$ on $M$ which is symmetric in $L^2(\mu)$, where $\mu:=
\e^{V(x)}\d x$ for $\d x$ the Riemannian volume measure.  Assume
that $\ll_0:= \inf\si(-L)$ is an eigenvalue of $-L$. Since $M$ is
connected, $\ll_0$ has a unique unit eigenfunction $\varphi_0>0$.

To clarify the  meaning in geometry analysis of the intrinsic
ultracontractivity, let us recall that

$$\|P_t^{\varphi_0}\|_{L^1(\mu_{\varphi_0})\to
L^\infty(\mu_{\varphi_0})}= \sup_{x,y\in M} \ff{ \e^{\ll_0 t}
h(x,y,t)}{\varphi_0(x)\varphi_0(y)},\ \ \ t>0,$$ where $h(x,y,t)$ is
the heat kernel of $P_t$ with respect to the weighted volume measure
$\mu.$

In order to study the intrinsic ultracontractivity of $P_t$, we make
use of the following intrinsic super-Poincar\'e inequality
introduced in \cite{W02} (see also \cite{OW}):

\beq\label{ISP} \mu(f^2) \le r \mu(|\nn f|^2) +\bb(r) \mu(\varphi_0
|f|)^2,\ \ \ r>0, f\in C_0^1(M),\end{equation} where $\bb:
(0,\infty)\to (0,\infty)$ is a decreasing function and

$$\mu(f):= \int_M f\d\mu,\ \ \ f\in L^1(\mu).$$

The intrinsic ultracontractivity of $P_t$ implies (\ref{ISP}) for
some $\bb$ (see \cite[Theorem 3.1]{W02}), and (\ref{ISP}) holds for
some $\bb$ if and only if $\si_{ess}(L)=\emptyset$ (see
\cite[Theorem 2.2]{W02}). On the other hand, if

\beq\label{1.2} \Psi(t):=\int_t^\infty \ff{\bb^{-1}(s)} s \d
s<\infty,\ \ \ t>\inf_{r>0}\bb(r),\end{equation} where
$\bb^{-1}(s):= \inf\{r>0: \bb(r)\le s\}$ for a positive decreasing
function $\bb$, then (\ref{ISP}) implies the intrinsic
ultracontractivity of $P_t$ with (see \cite[Theorem 3.3]{W00})

\beq\label{1.2'} \|P_t^{\varphi_0}\|_{L^1(\mu_{\varphi_0})\to
L^\infty(\mu_{\varphi_0})}\le \max\big\{\vv^{-1}\inf \bb,\
\Psi^{-1}((1-\vv)t)\big\}^2<\infty,\ \ \ \vv\in (0,1),
t>0.\end{equation} We refer to \cite{DS} for the study of intrinsic
ultracontractivity using the log-Sobolev inequality with parameters.

In section 2, (\ref{ISP}) with explicit $\bb$ is established in
terms of curvature lower bounds of $L$, the first Dirichlet
eigenvalue of $L$ outside large balls, and the super Poincar\'e
inequality on $M$ (see Theorem \ref{T2.1} below). As applications,
we obtain the following two concrete results (see Section 3 for
complete proofs).

To state these results, we introduce some curvature conditions. Let
$\Sec$ and $\Ric$ denote the sectional curvature and the Ricci
curvature on $M$ respectively. Let $\rr$ be the Riemannian distance
on $M$, and simply write $\rr_o:=\rr(o,\cdot)$ for a fixed reference
point $o\in M.$  Let $k$ and $K$ be two positive increasing
functions on $[0,\infty)$ such that

\beq\label{1.3} \Sec\le -k\circ \rr_o,\ \   \Ric\ge -K\circ\rr_o,\ \
\ \rr_o\gg 1
\end{equation} holds on $M$. Here $\Sec \le -k\circ\rr_o$ means that
for any $x\in M$ and any unit vectors $X,Y\in T_x$ with $\<X,Y\>=0$,
one has $\Sec(X,Y)\le -k(\rr_o(x));$ while $\Ric \ge -K\circ\rr_o$
means that $\Ric(X,X)\ge -K(\rr_o(x))|X|^2$ for any $x\in M$ and
$X\in T_x$. Finally, for a positive increasing function $h$ on
$(0,\infty)$, we let

$$h^{-1}(r):= \inf\{s>0:\ h(s)\ge r\},\ \ \ \ r\ge 0.$$

\beg{thm}\label{T1.1} Let $M$ be a Cartan-Hadamard manifold with
$d\ge 2$ and let $L=\DD$. Assume that $(\ref{1.3})$ holds for some
positive increasing functions $k$ and $K$ with $k(\infty)=\infty$.
We have: \beg{enumerate} \item[$(1)$] $(\ref{ISP})$ holds with

$$\bb(r):= \theta r^{-d/2} \exp\Big[ \theta k^{-1} (\theta/r)
\ss{K(4+ 2 k^{-1} (\theta/r))}\Big],\ \ \ r>0$$ for some constant
$\theta>0.$ \item[$(2)$] If

\beq\label{1.4} k^{-1}(R)\ss{K(4+2 k^{-1}(R))}\le c R^\vv,\ \ \ R\gg
1\end{equation} holds for some constants $c>0$ and $\vv\in (0,1)$,
then $P_t$ is intrinsically ultracontractive with

\beq\label{1.5} \|P_t^{\varphi_0}\|_{L^1(\mu_{\varphi_0})\to
L^\infty(\mu_{\varphi_0})}\le \exp\big[C(1+ t^{-\vv/(1-\vv)})\big],\
\ \ t>0\end{equation} for some constant $C>0,$ or equivalently

\beq\label{H} h(x,y,t)\le \e^{-\ll_0 t}
\varphi_0(x)\varphi_0(y)\exp\big[C(1+ t^{-\vv/(1-\vv)})\big],\ \
x,y\in M, t>0.\end{equation}

\item[$(3)$] If
 $(\ref{1.4})$ holds for some $c>0$ and $\vv=1$, then $P_t$ is
intrinsically hypercontractive.\end{enumerate}
\end{thm}

\beg{rem} (a) If $\Ric\ge -K$ for some constant $K\ge 0$, then
$\si_{ess}(\DD)\ne\emptyset.$ Since $M$ is non-compact and complete,
this follows from a comparison theorem by Cheng \cite{Cheng} for the
first Dirichlet eigenvalue and the Donnely-Li decomposition
principle \cite{DL} for the essential spectrum:

$$\inf \si_{ess}(-\DD)\le \sup_{x\in M} \ll_0(B(x,1))\le \ll_0(K),$$
where $\ll_0(B(x,1))$ is the first Dirichlet eigenvalue of $-\DD$ on
$D$ and $\ll_0(K)$ is the one on the unit geodesic ball in the
$d$-dimensional parabolic space with Ricci curvature equal to $K$.
Thus, the assumption $K(\infty)=\infty$ in Theorem \ref{T1.1} is
necessary for (\ref{ISP}) to hold. Correspondingly, the assumption
that $k(\infty)=\infty$ is also reasonable.

(b) The upper bound given in (\ref{1.5}), which is sharp due to
Example 1.1 below, is quite different from the known one on bounded
domains. Indeed, for $P_t$ the Dirichlet heat semigroup on a bounded
$C^{1,\aa} (\aa>0)$ domain in $\R^d$, the short time behavior of the
intrinsic  heat kernel is algebraic rather than exponential (see
\cite{OW}):

$$\sup_{x,y} \ff{h(x,y,t)\e^{\ll_0 t}}{\varphi_0(x)\varphi_0(y)}=
\bigcirc \big(t^{-(d+2)/2}\big).$$
\end{rem}

\

 Next, we consider the
case with drift. To this end, we adopt the following Bakry-Emery
curvature $\Ric_{m,L}$ instead of $\Ric$. Assume that for some $m>0$
and positive increasing  function $K$ one has, instead of the second
condition in (\ref{1.3}),

\beq\label{1.6} \Ric_{L,m}:= \Ric -\Hess_V- \ff {\nn V\otimes\nn V}m
\ge -K\circ \rr_o.\end{equation} Moreover, let $\ggg$ be a positive
increasing  function on $[0,\infty)$ such that

\beq\label{1.7} L\rr_o \ge \ss{\ggg\circ\rr_o},\ \ \ \rr_o\gg
1.\end{equation}

\beg{thm}\label{T1.2} Let $o$ be a pole in $M$ such that
$(\ref{1.6})$ and $(\ref{1.7})$ hold for some increasing positive
functions $K$ and $\ggg$ with $\ggg(\infty)=\infty$. Then
$\si_{ess}(L)=\emptyset$. Moreover, assuming

\beq\label{1.8} \lim_{\rr_o(x)\to\infty}
\ff{\ss{K(2+2\rr_o(x))}}{\log^+\mu(B(x,1))}=0,\end{equation} where
$B(x,1)$ is the unit geodesic ball at $x$, we have: \beg{enumerate}
\item[$(1)$]   $(\ref{ISP})$ holds
with

$$\bb(r)= \theta r^{-(m+d+1)/2} \exp\Big[\theta \ggg^{-1}
(32/r)\ss{K(2 +2 \ggg^{-1} (32/r))}\Big], \ \ \ r>0$$ for some
constant $\theta >0$. \item[$(2)$] If there exist  $c>0$ and $\vv\in
(0,1)$ such that

\beq\label{1.9}\ggg^{-1}(R)\ss{K(2+2 \ggg^{-1}(R))}\le c R^\vv,\ \ \
R\gg 1,\end{equation} then $P_t$ is intrinsically ultracontractive
with $(\ref{1.5})$ and $(\ref{H})$  holding for some constant $C>0.$
\item[$(3)$] If
 $(\ref{1.9})$ holds for some $c>0$ and $\vv=1$, then $P_t$ is
intrinsically hypercontractive.\end{enumerate}
\end{thm}

To conclude this section, we present below two typical examples to
show that conditions in Theorems \ref{T1.1} and \ref{T1.2} are
sharp. To make the introduction brief, we leave their proofs to
Section 4.

\paragraph{Example 1.1.} Let $M$ be a Cartan-Hadamard manifold with

$$ -c_1 \rr_o^\dd\le \Sec\le -c_2 \rr_o^\dd,\ \ \ \rr_o\gg 1$$ for some constants
$c_1, c_2,\dd >0$. Then $\si_{ess}(\DD)=\emptyset$ and for $L=\DD$,
(\ref{ISP}) holds with

$$\bb(r)= \exp[c(1+r^{-(2+\dd)/[2\dd]})]$$ for some constant $c>0.$
Consequently:

\begin{enumerate}
\item[(1)] $P_t$ is intrinsically ultracontractive if and only if
$\dd>2$, and when $\dd>2$ one has

$$\|P_t^{\varphi_0}\|_{L^1(\mu_{\varphi_0})\to
L^\infty(\mu_{\varphi_0})}\le \theta_1 \exp\big[\theta_2
t^{-(\dd+2)/(\dd-2)}\big],\ \ \ t>0$$ for some constants
$\theta_1,\theta_2>0$, which is sharp in the sense that the constant
$\theta_2$ can not be replaced by  any positive function
$\theta_2(t)$ with $\theta_2(t)\downarrow 0$ as $t\downarrow 0$.

\item[(2)] $P_t$  is intrinsically hypercontractive if and only if
$\dd\ge 2$.
\end{enumerate}

\paragraph{Example 1.2.} Let $M$ be a Cartan-Hadamard manifold with

$$ \Ric\ge -c (\rr_o^{2(\dd-1)}+1)$$ for some constants
$c>0$ and $\dd >1$. Let $V= \theta \rr_o^\dd$ for some constant
$\theta
>0$ and $\rr_o\gg 1.$ Then $\si_{ess}(L)=\emptyset$ and (\ref{ISP}) holds with

$$\bb(r)= \exp[c(1+r^{-\dd/[2(\dd-1)]})]$$ for some constant $c>0.$
Consequently:

\begin{enumerate}
\item[(1)] $P_t$ is intrinsically ultracontractive if and only if
$\dd>2$, and when $\dd>2$ one has

$$\|P_t^{\varphi_0}\|_{L^1(\mu_{\varphi_0})\to
L^\infty(\mu_{\varphi_0})}\le \theta_1 \exp\big[\theta_2
t^{-\dd/(\dd-2)}\big],\ \ \ t>0$$ for some constants
$\theta_1,\theta_2>0$, which is sharp in the sense that the constant
$\theta_2$ can not be replaced by  any positive function
$\theta_2(t)$ with $\theta_2(t)\downarrow 0$ as $t\downarrow 0$.
\item[(2)] $P_t$  is intrinsically hypercontractive if and only if
$\dd\ge 2$.
\end{enumerate}

\section{The intrinsic super-Poincar\'e inequality}

As explained in the last section that due to \cite[Theorem
2.2]{W02}, (\ref{ISP}) holds for some $\bb$ if and only if
$\si_{ess}(L)=\emptyset$. According to the Donnely-Li decomposition
principle (see \cite{DL}), they are also equivalent to

$$\ll_0(R):= \inf\{\mu(|\nn f|^2):\ \mu(f^2)=1, f\in C_0^1(M),
f|_{B(o,R)}=0\} \uparrow \infty$$ as $R\uparrow\infty.$ The purpose
of this section is to estimate $\bb$ in (\ref{ISP}) by using
$\ll_0(R)$ and the curvature condition. To this end, we will make
use of the following super-Poincar\'e inequality:

\beq\label{S} \mu(f^2)\le r \mu(|\nn f|^2) + \bb_0(r)\mu(|f|)^2,\ \
\ r>0, f\in C_0^1(M)\end{equation} for some decreasing function
$\bb_0: (0,\infty)\to (0,\infty).$ In particular, by \cite[Corollary
1.1 (2)]{W00}, (\ref{S}) with $\bb_0(r)= c(1+r^{-p/2})$ for some
constant $c>0$ and $p>2$ is equivalent to the classical Sobolev
inequality

$$\mu(|f|^{2p/(p-2)})^{(p-2)/p}\le C(\mu(|\nn f|^2) +\mu(f^2)),\ \ \
f\in C_0^1(M)$$ for some constant $C>0.$ The latter inequality holds
for a large class of non-compact manifolds. For instance, according
to \cite{Croke}, it holds true for $V=0$ provided either the
injectivity radius of $M$ is infinite, or the injectivity radius is
positive and the Ricci curvature is bounded below.

\beg{thm} \label{T2.1} Assume $(\ref{S})$. Let $K$ be positive
increasing function on $[0,\infty)$ such that $(\ref{1.6})$ holds.
If $\ll_0(R)\uparrow\infty$ as $R\uparrow\infty$, then $(\ref{ISP})$
holds with

$$\bb(r)=C \bb_0(r/8)\exp\Big[ C \ll_0^{-1}(8/r) K(2+2\ll_0^{-1}
(8/r))\Big],\ \ \ r>0.$$\end{thm} To prove this result, we first
estimate the ground state $\varphi_0$ from below. The following
lemma is proved by using the Li-Yau (\cite{LY}) type parabolic
Harnack inequality derived by X.-D. Li \cite{Li}.

\beg{lem} \label{L2.2} If $(\ref{1.6})$ holds then for the positive
ground state $\varphi_0$, there exists a constant $C>0$ such that

$$\varphi_0\ge \ff 1 C \exp\Big[-C\rr_o
\ss{K(2\rr_o)}\Big].$$\end{lem}

\beg{proof} Since $\varphi_0$ is bounded below by a positive
constant on a compact set, it suffices to prove for $\rr_o\ge 1.$
Let $x\in M$ with $\rr_o(x)\ge 1.$ Applying \cite[Theorem 5.2]{Li}
to $\aa=2$ and $R=\rr_o(x)$, we obtain

\beg{equation*}\beg{split} \e^{-\ll_0} \varphi_0(o) = P_1
\varphi_0(o) &\le (P_{1+s}\varphi_0(x))
(1+s)^{m+d}\exp\Big[c_1^2sK(2\rr_o(x))+ \ff{\rr_o(x)^2}{2s}\Big]\\
&= \varphi_0(x) \e^{-\ll_0(1+s)}(1+s)^{m+d}
\exp\Big[c_1^2sK(2\rr_o(x))+ \ff{\rr_o(x)^2}{2s}\Big],\ \ \
s>0\end{split}\end{equation*} for some constant $c_1>0.$ Then the
proof is completed by taking $s= \rr_o(x)/\ss{K(2\rr_o(x))}.$
\end{proof}

\ \newline \emph{Proof of Theorem \ref{T2.1}.} Since one may always
take decreasing $\bb$, it suffices to prove for $r\le 1.$ Let $f\in
C_0^1(M)$ be fixed. Let $h_R= (\rr_o-R)^+\land 1,\ R>0.$ Then $h_R$
is Lipschitz continuous so that (\ref{S}) applies to $f(1-h_R)$
instead of $f$:

\beq\label{2.1} \mu(f^2(1-h_R)^2)\le 2 s\mu(|\nn f|^2) + 2 s\mu(f^2)
+\bb_0(s) \mu(|f|1_{B(o, R+1)})^2,\ \ \ s>0.\end{equation} Next,
since $h_Rf=0$ on $B(o, R)$, we have

$$\mu(f^2h_R^2)\le \ff {\mu(|\nn (fh_R)|^2)} {\ll_0(R)} \le \ff 2
{\ll_0(R)} \mu(|\nn f|^2) + \ff 2 {\ll_0(R)} \mu(f^2).$$ Combining
this with (\ref{2.1}) we obtain

\beg{equation*}\beg{split} \mu(f^2) &\le 2 \mu(f^2 h_R^2)+ 2
\mu(f^2(1-h_R)^2)\\
&\le \Big(4s +\ff 4 {\ll_0(R)}\Big)(\mu(|\nn f|^2)+\mu(f^2))
+\ff{2\bb_0(s)}{\inf_{B(o,R+1)}\varphi_0^2}\mu(|f|\varphi_0)^2.\end{split}\end{equation*}
Thus, if $4s +\ff 4 {\ll_0(R)}\le \ff 1 2$ then

$$\mu(f^2) \le \Big(8s +\ff 8 {\ll_0(R)}\Big) \mu(|\nn f|^2)
+\ff{4\bb_0(s)}{\inf_{B(o,R+1)}\varphi_0^2}\mu(|f|\varphi_0)^2.$$
Hence, (\ref{ISP}) holds for

$$\bb(r):= \inf \bigg\{\ff{4\bb_0(s)}{\inf_{B(o,R+1)}\varphi_0^2}:\ 8s +\ff 8 {\ll_0(R)}\le r\bigg\},\ \ \
r\le 1.$$ Combining this with Lemma \ref{L2.2}, there exists a
constant $c>0$ such that (\ref{ISP}) holds for

$$\bb(r):= \inf \bigg\{c\bb_0(s)\exp\big[c(R+1)\ss{K(2+2R)}\big]:\ 8s +\ff 8 {\ll_0(R)}\le r\bigg\},\ \ \
r\le 1.$$ This completes the proof by taken $s= r/8$ and $R=
\ll_0^{-1} (8/r).$ \qed

\section{Proofs of Theorems \ref{T1.1} and \ref{T1.2}}

\emph{Proof of Theorem \ref{T1.1}.} a) Since $V=0$, (\ref{1.3})
implies (\ref{1.6}). Moreover, since $M$ is a Cartan-Hadamard
manifold, its injectivity radius is infinite. Hence, by
\cite{Croke}, one has $\|P_t\|_{L^1(\mu)\to L^\infty(\mu)}\le c
t^{-d/2}$ for some $c>0$ and all $t>0.$ By \cite[Theorem 4.5
(b)]{W00}, this implies (\ref{S}) with $\bb_0(r)= c(1+r^{-d/2})$ for
some constant $c>0.$

b) Since $M$ is a Cartan-Hadamard manifold, $B(o,R)^c$ is concave.
Let $R_0>0$ be such that   (\ref{1.3}) holds for $\rr_o\ge R_0$.
Then for any $R\ge R_0$, we have  $\Sec\le -k(R)$ on $B(o,R)^c.$ To
make use of the Laplacian comparison theorem, we note that the
distance to the boundary of $B(o,R)^c$ is $\rr_o-R$ for $\rr_o\ge
R$, and that the function

$$h(s):=\cosh\big(\ss{k(R)}\, s\big),\ \ \ s\ge 0$$ solves the equation

$$h''(s) -k(R)h(s) =0,\ \ h(0)=1,\ h'(0)=0.$$ Then  (see \cite[Theorem
0.3]{K})

\beq\label{3.1} \DD\rr_o \ge \ff{(d-1)h'(\rr_o-R)}{h(\rr_o-R)}\ge
c_0 \ss{k(R)},\ \ \ \rr_o\ge R+1\end{equation} holds for some
constant $c_0>0$ which is independent of $R$. By the Green formula,
for any smooth domain $D\subset B(o,R+1)^c$, it follows from
(\ref{3.1}) that

$$c_0\ss{k(R)}\mu(D)\le \int_D \DD \rr_o\d\mu \le \int_{\pp D}
|N\rr_o | \d\mu_\pp\le \mu_\pp(\pp D),$$ where $N$ is the unit
normal vector field on $\pp D$ and $\mu_\pp$ is the  measure on $\pp
D$ induced by $\mu$. Thus, by Cheeger's inequality (see
\cite{Cheeger}), we have

$$\ll_0(R+1)\ge \ff {c_0^2 k(R)}4,\ \ \ R\ge R_0$$ which goes to infinite as $R\to\infty$, so that
$\si_{ess}(\DD)=\emptyset.$  Moreover,

$$\ll_0^{-1}(8/r)\le \inf\Big\{R+1: R\ge R_0, \ff {c_0^2} 4 k(R)\ge \ff 8
r\Big\} = 1 +R_0\lor k^{-1} (32/c_0^2r),\ \ \ r>0.$$ Then by Theorem
\ref{T2.1} with $\bb_0(r)=
  c(1+r^{-d/2})$, we obtain the desired $\bb(r)$ for  some
  $\theta>0$.

c) If (\ref{1.4}) holds then by (1),  (\ref{ISP}) holds for

$$\bb(r)= \exp[\theta(1+ r^{-\vv})]$$ for some constant $\theta>0.$
If $\vv\in (0,1)$ then (\ref{1.5}) follows from \cite[Corollary
3.4(1)]{W02}.  If $\vv=1$ then

$$\mu(f^2)\le r \mu(|\nn f|^2) +\exp[\theta(1+ r^{-1})]
\mu(\varphi_0 |f|)^2,\ \ \ r>0, f\in C_0^1(M).$$ Applying this to
$f\varphi_0$ and noting that

$$ \ff 1 2 \mu(\<\nn f^2, \nn \varphi_0^2\>)= -\ff 1 2
\mu(f^2L\varphi_0^2)= \ll_0  \mu_{\varphi_0}(f^2)- \mu(f^2|\nn
\varphi_0|^2),$$
 we arrive at

\beg{equation*}\beg{split} \mu_{\varphi_0}(f^2) &\le r
\mu_{\varphi_0}(|\nn f|^2) + r \mu(f^2|\nn \varphi_0|^2) +\ff r 2
\mu(\<\nn f^2, \nn \varphi_0^2\>) +\e^{\theta(1+r^{-1})}
\mu_{\varphi_0}(|f|)^2\\
&= r\mu_{\varphi_0}(|\nn f|^2) +  r\ll_0  \mu_{\varphi_0}(f^2)
+\e^{\theta(1+r^{-1})} \mu_{\varphi_0}(|f|)^2,\ \ \
r>0.\end{split}\end{equation*} This implies

$$ \mu_{\varphi_0}(f^2) \le 2 r\mu_{\varphi_0}(|\nn f|^2) +2\e^{\theta(1+r^{-1})}
\mu_{\varphi_0}(|f|)^2,\ \ \ r\in (0, 1/(2\ll_0)).$$ Hence, there
exists a constant $\theta'>0$ such that

\beq\label{SS}\mu_{\varphi_0}(f^2) \le  r\mu_{\varphi_0}(|\nn f|^2)
+\e^{\theta'(1+r^{-1})} \mu_{\varphi_0}(|f|)^2,\ \ \
r>0.\end{equation} By \cite[Corollary 1.1(1)]{W00}, this is
equivalent to the defective log-Sobolev inequality

\beq\label{D} \mu_{\varphi_0}(f^2\log f^2)\le
C_1\mu_{\varphi_0}(|\nn f|^2) + C_2,\ \ \ f\in C_b^1(M),
\mu_{\varphi_0}(f^2) =1 \end{equation}  for some $C_1,C_2>0.$ On the
other hand, (\ref{SS}) and the weak Poincar\'e inequality due to
\cite[Theorem 3.1]{RW} imply the Poincar\'e inequality (see
\cite[Proposition 3.1]{RW})

$$\mu_{\varphi_0}(f^2) \le C\mu_{\varphi_0}(|\nn f|^2)
+ \mu_{\varphi_0}( f )^2,\ \ \ f\in C_b^1(M)$$ for some constant
$C>0$. Combining this and (\ref{D}) we obtain the strict log-Sobolev
inequality, namely, (\ref{D}) with $C_2=0$ and   some possibly
different $C_1>0$. Therefore, due to \cite{G}, $P_t^{\varphi_0}$ is
hypercontracive since it is associated to the Dirichlet form
$\mu_{\varphi_0}(\<\nn \cdot,\nn \cdot\>)$ on
$H^{2,1}(\mu_{\varphi_0})$.  We remark that the implication of the
hypercontactivity from the defective log-Sobolev inequality can also
be deduced by using the uniformly positivity improving property of
the diffusion semigroup, see e.g. \cite{A} for details. Then the
proof is finished.\qed

\

 To prove Theorem \ref{T1.2}, we first establish the super
Poincar\'e inequality (\ref{S}) for a concrete $\bb_0$.

\beg{lem}\label{L3.1} In the situation of Theorem $\ref{T1.2}.$
$(\ref{1.8})$ implies $(\ref{S})$ with $\bb_0(r)=
c(1+r^{-(m+d+1)/2})$ for some constant $c>0.$ \end{lem}

\beg{proof} By \cite[Theorem 5.2]{Li} with $\aa= (m+d+1)/(m+d)$, for
any measurable function $f\ge 0$ with $\mu(f)=1,$ we have

$$P_t f(x)\le (P_{t+s}f(y))\Big(1+ \ff s
t\Big)^{(m+d+1)/2}\exp\Big[CK(2[\rr_o(x)\lor\rr_o(y)])s +\ff {\aa
\rr(x,y)^2}{4s}\Big]$$ for some constant $C>0$ and all $s,t>0.$ This
implies

$$1 = \int_M P_{t+s}f(y)\mu(\d y) \ge (P_t f(x))\Big(1+ \ff s
t\Big)^{-(m+d+1)/2}\int_{B(x,1)}\e^{-CK(2[\rr_o(x)\lor\rr_o(y)])s -
\aa /[4s]}\mu(\d y).$$ Taking $s= 1/\ss{K(2+2\rr_o(x))},$ we obtain

$$P_t f(x)\le c_0(1+t^{-1})^{(m+d+1)/2} \cdot\ff{\exp[c_0
\ss{K(2+2\rr_o(x))}]}{\mu(B(x,1))}$$ for some constant $c_0>0$ and
all $t>0, x\in M.$  Combining this with (\ref{1.8}) we obtain

$$\|P_t\|_{L^1(\mu)\to L^\infty(\mu)}\le c_1(1+
t^{-1})^{(m+d+1)/2},\ \ \ t>0$$ for some constant $c_1>0.$ According
to \cite[Theorem 4.5(b)]{W00}, this is equivalent to (\ref{S}) with
$\bb_0(r)=c(1+r^{-(m+d+1)/2})$ for some constant $c>0.$ \end{proof}

\ \newline \emph{Proof of Theorem \ref{T1.2}.} By (\ref{1.7}) and
Cheeger's inequality explained in b) of the proof of Theorem
\ref{T1.1}, we have

$$\ll_0(R)\ge \ff {\ggg(R)}4,\ \ \ R\gg 1.$$ Since $\ggg(R)\to\infty$
as $R\to\infty$, the essential spectrum of $L$ is empty and the
desired $\bb$ follows from Theorem \ref{T2.1} and Lemma \ref{L3.1}.
The remainder of the proof is then similar to that of Theorem
\ref{T1.1}.\qed

\section{Proofs of Examples 1.1 and 1.2}

\emph{Proof of Example 1.1.} Since $\Sec\le -c_2 \rr_o^\dd$ for some
$c_2,\dd>0$ and large $\rr_o$, Theorem \ref{T1.1} implies
$\si_{ess}(\DD)=\emptyset.$ Moreover, one may take $K(r)= (d-1)c_1
r^\dd$ and $k(r)= c_2 r^\dd$ for large $r$, so that

$$k^{-1}(R)\ss{K(4+2 k^{-1}(R))}\le c R^{\ff 1 2 +\ff 1 \dd}$$ for
some constant $c>0$ and large $R$. Then the sufficiency and the
desired upper bound of $\|P_t^{\varphi_0}\|_{L^1(\mu_{\varphi_0})\to
L^\infty(\mu_{\varphi_0})}$ follow from Theorem \ref{T1.1}.

Next, by the concrete $K$ and Lemma \ref{L2.2} we have

\beq\label{4.1} \varphi_0 \ge \ff 1 C
\exp\big[-C\rr_o^{1+\dd/2}\big]\end{equation} for some constant
$C>0.$ If $P_t$ is intrinsically ultracontractive, i.e.
$P_t^{\varphi_0}$ is ultracontractive by definition, then, according
to \cite[Theorem 2.2.4]{Davies} (see also \cite{G} and \cite{DS}),
there exists a function $\bb: (0,\infty)\to (0,\infty)$ such that

$$\mu_{\varphi_0}(f^2\log f^2)\le r\mu_{\varphi_0}(|\nn f|^2) + \bb(r),\ \ \
f\in C_b^1(M), \mu_{\varphi_0}(f^2)=1.$$ By the concentration of
reference measures induced by log-Sobolev inequalities (see e.g.
  \cite[Corollary 6.3]{RW03}), the above log-Sobloev inequality implies  $\mu_{\varphi_0}(\e^{\ll
\rr_o^2})<\infty$ for any $\ll>0.$ Combining this with (\ref{4.1})
and noting that the Riemannian volume of a Cartan-Hadamard manifold
is infinite, we conclude that $\dd>2.$ Similarly, if $P_t$ is
intrinsically hypercontractive, then $\mu_{\varphi_0}(\e^{\ll
\rr_o^2})<\infty$ for some $\ll>0,$ so that $\dd\ge 2.$

Finally, let $\dd>2.$ If there exists $\theta_1>0$ and a positive
function $h$ with $h(t)\downarrow 0$ as $t\downarrow 0$ such that
$$\|P_t^{\varphi_0}\|_{L^1(\mu_{\varphi_0})\to
L^\infty(\mu_{\varphi_0})}\le \theta_1 \exp\big[\theta_2
t^{-(\dd+2)/(\dd-2)}\big],\ \ \ t>0,$$ then  \cite[Theorem 4.5]{W00}
implies (\ref{ISP})  for

\beg{equation*}\beg{split} \bb(r)&= \inf_{s\le r, t>0} \ff{s} t
\|P_t^{\varphi_0}\|_{L^1(\mu_{\varphi_0})\to
L^\infty(\mu_{\varphi_0})} \e^{t/s -1}\\
&=\theta_1 \inf_{s\le r, t>0} \ff s t
\exp\Big[h(t)t^{-(\dd+2)/(\dd-2)}+\ff t s -1\Big],\ \ \
r>0.\end{split}\end{equation*} Taking $s=r\land 1$ and $t=
\big(r^{(\dd-2)/(2\dd)}\land 1\big) h(r^{(\dd-2)/(2\dd)}\land
1)^{(\dd-2)/(2\dd)},$ we obtain

\beq\label{FF} \bb(r)\le \theta_2 \exp\big[\tt h(r)
r^{-(\dd+2)/(2\dd)}\big],\ \ \ r>0\end{equation} for some constant
$\theta_2>0$ and positive function $\tt h$ with $\tt h(r)\downarrow
0$ as $r\downarrow 0$.

Finally, we aim to deduce from (\ref{FF})  that

\beq\label{FF2} \mu(\e^{\ll \rr_o^{1+\dd /2}})<\infty,\ \ \
\ll>0,\end{equation} which is contradictive to (\ref{4.1}). To this
end, we apply \cite[Theorem 6.2]{W00b}, which says that

\beq\label{CC}\mu(\exp[c_1 \rr_o\xi(c_2
\rr_o)])<\infty\end{equation} holds for some constants $c_1,c_2>0$
and

$$\xi(\ll):= \inf\bigg\{s\ge 1:\ \int_1^s\ff 1 {t^2}\log \bb(1/(2 t^2))\d
t \ge \ll\bigg\},\ \ \ \ll>0.$$ Since  (\ref{FF}) implies

\beg{equation*}\beg{split} \int_1^s \ff 1 {t^2}\log \bb(1/(2 t^2))\d
t&\le \theta_3 +\theta_3 \int_1^s t^{-(\dd-2)/\dd} \tt h(1/(2
t^{(\dd-2)/\dd}))\d t\\
&\le \theta_3 + \vv(s) s^{2/\dd},\ \ \ \
s>1\end{split}\end{equation*}  for some constant $\theta_3>0$ and
some positive function $\vv$ with $\vv(s)\downarrow 0$ as $s\uparrow
\infty,$ one has $\xi(\ll) \ll^{-\dd/2}\to\infty$ as $\ll\to\infty$.
Therefore, (\ref{FF2}) follows from (\ref{CC}).\qed

\ \newline \emph{Proof of Example 1.2.}  Since $M$ is a
Cartan-Hadamard manifold and $\dd>1$, we have

$$L\rr_o\ge \dd \rr_o^{\dd-1}=:\ss{\ggg\circ\rr_o},\ \ \ \rr_o\gg 1.$$
In particular, $\ggg(\infty)=\infty$ so that
$\si_{ess}(L)=\emptyset.$ Moreover, since

$$\Ric\ge -c(1+\rr_o^{2(\dd-1)}),\ \ \ |\nn V|^2= \theta^2 \dd^2
\rr _o^{2(\dd-1)}$$ and $\Hess_V= \theta \Hess_{\rr_o^\dd}\ge 0$ for
large $\rr_o$ as $M$ is Cartan-Hadamard, we may take $K(r)=
c_1(1+r^{2(\dd-1)})$ for some constant $c_1>0.$ Therefore,
(\ref{1.9}) holds for some $c>0$ and $\vv= \ff 1 2 +\ff 1
{2(\dd-1)}.$  Therefore, the sufficiency follows from Theorem
\ref{T1.2} provided (\ref{1.8}) holds. Indeed, since $M$ is a
Cartan-Hadamard manifold, we have

$$\mu(B(x,1))\ge c(d) \exp\Big[\inf_{B(x,1)}V\Big]\ge c(d)
\exp[\theta (\rr_o(x)-1)^\dd],\ \ \ \rr_o(x)\ge 1,$$ where $c(d)$ is
the volume of the unit ball in $\R^d$. This implies (\ref{1.8}).

On the other hand, by Lemma \ref{L2.2} and the concrete $K$, we have

$$\varphi_0 \ge \ff 1 C \exp\big[-C\rr_o^\dd\big]$$ for some constant
$C>0.$ Then the remainder of the proof is as same as that in the
proof of Example 1.1. \qed

\paragraph{Acknowledgement.} The author  would like to thank the
referees   for   useful comments.

 \beg{thebibliography}{99}

\bibitem{A} S. Aida, \emph{Uniformly positivity improving property,
Sobolev inequalities and spectral gap,} J. Funct. Anal. 158 (1998),
152--185.

\bibitem{Cheeger} J. Cheeger,  \emph{A lower bound for the smallest
eigenvalue of the Laplacian,} Problems in analysis, a symposium in
honor of S. Bochner 195--199, 1970: Princeton U. Press, Princeton.

\bibitem{CS} Z.-Q. Chen and R. Song, \emph{Intrinsic
ultracontractivity, conditional lifetimes and conditional gauge
for symmetric stable processes on rough domains,} Illinois J.
Math. 44 (2000), 138-160.

\bibitem{Cheng} S. Y. Cheng, \emph{Eigenvalue comparison theorems
and its geometric applications,} Math. Z. 143(1975), 289--297.

\bibitem{C} F. Ciprina, \emph{Intrinsic ultracontractivity of Dirichlet Laplacians in nonsmooth domains,}
Potential Anal. 3 (1994), 203-218.

\bibitem{Croke} C. B. Croke, \emph{Some isoperimetric inequalities and eigenvalue
estimates,} Ann. Sci. \'Ec. Norm. Super. 13(1980), 419--435.

\bibitem{Davies} E. B. Davies, \emph{Heat Kernels and Spectral
Theory,} Cambridge University Press, Cambridge, 1989.

\bibitem{DS}  E. B. Davies and B. Simon, \emph{Ultracontractivity and heat kernel
for Schr\"odinger operators and Dirichlet Laplacians,} J. Funct.
Anal. 59(1984), 335--395.

\bibitem{DL} H. Donnelly and P. Li, \emph{Pure point spectrum and negative
curvature for noncompact manifolds,} Duke Math. Journal 46(1979),
497--503.

\bibitem{G} L. Gross,  \emph{Logarithmic Sobolev inequalities,} Amer. J. Math.
{\bf 97}(1976), 1061--1083.

\bibitem{GR} L. Gross and O. Rothaus, \emph{Herbst inequalities
for supercontractive semigroups,} J. Math. Kyoto Univ. 38(1998),
295--318.

\bibitem{K} A. Kasue, \emph{Applications of Laplacian and
Hessian comparison theorems,} Adv. Stud. Pure Math. 3(1984)
333--386.

\bibitem{KS} P. Kim and R. Song, \emph{Intrinsic ultracontractivity
for non-symmetric L\'evy processes,} preprint.

\bibitem{LY} P. Li and  S.-T. Yau, \emph{On the parabolic kernel of
the Schr\"odinger operator,} Acta Math. 156(1986), 153--201.

\bibitem{Li} X.-D. Li, \emph{Liouville theorems for symmetric
diffusion operators on complete Riemannian manifolds,} J. Math.
Pures Appl. 84(2005), 1295--1361.

\bibitem{OW} E. M. Ouhabaz and F.-Y. Wang, \emph{Sharp estimates for  intrinsic ultracontractivity on
$C^{1,\alpha}$-domains,} to appear in Manuscript Math.

\bibitem{RW} M. R\"ockner, M. and F.-Y. Wang,   \emph{Weak Poincar\'e inequalities
and $L^2$-convergence rates of Markov semigroups,} J. Funct. Anal.
185(2001), 564--603.

\bibitem{RW03}  M. R\"ockner and F.-Y. Wang, \emph{Supercontractivity and ultracontractivity for
(non-symmetric) diffusion semigroups on manifolds,} Forum Math.
15(2003), 893--921.

\bibitem{W00b} F.-Y. Wang, \emph{Functional inequalities for empty essential
spectrum,} J. Funct. Anal. 170(2000), 219--245.

\bibitem{W00} F.-Y. Wang, \emph{Functional inequalities, semigroup properties
and spectrum estimates,} Infin. Dimens. Anal. Quant. Probab.
Relat. Topics 3(2000), 263--295.

\bibitem{W02} F.-Y. Wang, \emph{Functional inequalities and spectrum estimates:
the infinite measure case,} J. Funct. Anal. 194(2002), 288--310.

\end{thebibliography}

 \end{document}